\documentclass[conference]{IEEEtran}
\IEEEoverridecommandlockouts
\usepackage{cite}
\usepackage{amsmath,amssymb,amsfonts}
\usepackage{algorithm}    

\usepackage{algpseudocode}

\usepackage{graphicx}
\usepackage{textcomp}
\usepackage{xcolor}
\usepackage{epsfig}
\usepackage{subfig}
\usepackage{algorithm}        
\usepackage{algpseudocode}    

\graphicspath{{pictures/}}
\def\BibTeX{{\rm B\kern-.05em{\sc i\kern-.025em b}\kern-.08em
    T\kern-.1667em\lower.7ex\hbox{E}\kern-.125emX}}
\begin{document}
\title{Soft Actor-Critic with Backstepping-Pretrained DeepONet for control of  PDEs
\thanks{The paper is supported by the National Natural Science Foundation of China (62173084), the Project of Science and Technology Commission of Shanghai Municipality, China (23ZR1401800).}
\thanks{Corresponding author: Jie Qi (jieqi@dhu.edu.cn)}
}

\author{
    \IEEEauthorblockN{1\textsuperscript{st} Chenchen Wang}
    \IEEEauthorblockA{\textit{College of Information Science}\\ \textit{ and Technology Donghua University} \\ Shanghai, China \\ 2242131@mail.dhu.edu.cn}
    \and
    \IEEEauthorblockN{2\textsuperscript{nd} Jie Qi}
    \IEEEauthorblockA{\textit{College of Information Science}\\ \textit{ and Technology Donghua University} \\ Shanghai, China \\ jieqi@dhu.edu.cn}
    \and
    \IEEEauthorblockN{3\textsuperscript{rd} Jiaqi Hu}
    \IEEEauthorblockA{\textit{College of Information Science}\\  \textit{ and Technology Donghua University} \\ Shanghai, China \\ hujiaqienjoy@gmail.com}
}

\maketitle

\begin{abstract}
This paper develops a reinforcement learning-based controller for the stabilization of partial differential equation (PDE) systems. Within the Soft Actor-Critic (SAC) framework, we embed a DeepONet, a well-known neural operator (NO), which is pretrained using the backstepping controller. The pretrained DeepONet captures the essential features of the backstepping controller and is used to extract features, replacing the convolutional neural networks (CNNs) in the original actor-critic networks, and directly connects to the fully connected layers of the SAC architecture. We apply this novel backstepping and reinforcement learning integrated method to stabilize an unstable 1D hyperbolic PDE and an unstable reaction-diffusion PDE. The proposed method is shown to outperform standard SAC in simulations, SAC with an untrained DeepONet, and the backstepping controller on both systems.

\end{abstract}


\section{Introduction}

Control of PDE systems remains a highly challenging task due to the infinite-dimensional nature of the state space and the complexity of system dynamics. Learning-based approaches integrated with the classical control theory have emerged as a promising alternative, offering data-driven adaptability \cite{quartz2024stochastic}. 

For finite-dimensional systems, there has been a growing body of research \cite{rodriguez2022lyanet,ip2024lyapunov,berkenkamp2017safe}. For instance, a control-informed reinforcement learning (RL) framework that integrates PID control components into the architecture of deep RL policies is proposed in \cite{Bloor2025}, 
In this \cite{yang2023certifiably}, it embeds barrier function theory into neural ordinary differential equations (ODEs), ensuring that the system state remains within a safe region by optimizing the parameters of a learning-based barrier function.

The infinite-dimensional nature of PDE systems poses significant challenges for control design, which has limited research progress in this field. For some special PDE systems there are some RL-based control design, such as applying proximal policy optimization to optimize the boundary controller for the Aw-Rascle-Zhang (ARZ) model of highway traffic flow in \cite{yu2021reinforcement}. Furthermore, by combining the state-dependent Riccati equation (SDRE) with Bayesian linear regression, \cite{alla2024online} proposes a RL-based algorithm for online identification and stable control of unknown PDEs.
 

However, one of the key challenges in combining classical control theory with learning-based methods lies in effectively incorporating prior knowledge from classical control into neural network learning.
To address this issue, we propose a novel approach where a NO, specifically a Deep Operator Network (DeepONet) \cite{lu2021learning}, \cite{71hu2024boundary}, is pretrained to learn the backstepping controller, a well-established classical control law for PDE systems. DeepONet is particularly suited for this task as it can represent mappings from function spaces to function spaces, making it ideal for approximating the feedback operator. The use of DeepONet to learn backstepping gain functions/controllers is first proposed by \cite{bhan2023neural,krstic2024neural}. Subsequently, this approach is extended to handle PDE systems with delay \cite{qi2024neural,qi2024neural2,wang2025deep}, applied to enhance the real-time performance of adaptive control of PDEs \cite{Lamarque2025,bhanneural}, and utilized in traffic flow control problems \cite{yuh2025neural}.


In our method, the DeepONet pretrained on the backstepping controller is used for feature extraction, as a replacement for the CNNs in standard SAC. The DeepONet connects directly to the fully connected actor-critic networks layers. Throughout the training phase, parameters of the DeepONet feature extraction module are optimized in conjunction with those of the SAC actor-critic networks, thereby enabling fine-tuning during reinforcement learning.

The primary contribution of this paper is the introduction of a backstepping-pretrained DeepONet into the RL architecture. This accelerates training by providing better initialization, allowing the reward function to start from a higher baseline. We apply this novel backstepping and RL integrated method to stabilize an unstable first-order hyperbolic PDE and an unstable reaction-diffusion PDE. Simulation results show that the proposed method significantly reduces transient oscillations compared to the backstepping controller, vanilla SAC, and SAC with an untrained DeepONet. It also decreases the steady-state error relative to both SAC baseline algorithms. In our experiments, RL-based controllers exhibit small steady-state errors, likely due to the stochastic nature of their policies, whereas the rigorous backstepping controller can eliminate such errors. By integrating a pretrained DeepONet infused with backstepping knowledge, the proposed method markedly reduces this error.

Additionally, when training the DeepONet, we incorporate not only state variables but also system coefficient functions as inputs. This design enables the pretrained DeepONet to adapt its output, the control effort, to variations in system coefficients. As a result, the trained RL controller remains effective even when applied to PDE systems with coefficients different from those seen during training, demonstrating the robustness of the method to parameter variations.

The structure of this paper is as follows: Section II provides problem formulation; Section III presents the NO-based RL framework; Section IV validates the effectiveness of the method through simulations; and finally, Section V concludes the paper.
\section{Problem description}
\subsection{General PDE control problem}
We consider a one-dimensional PDE defined on a domain  $\mathcal{X} = [0, 1]$. The time domain is $ \mathcal{T} = [0, T] \subset \mathbb{R}^+ $. Let $ u(x, t) $, where $ x \in \mathcal{X} $ and $ t \in \mathcal{T} $, describe the state of the system governed by the PDE, with the dynamics given by
\begin{align} 
	\label{PDE}
	\frac{\partial u}{\partial t} = F\left(u, \frac{\partial u}{\partial x}, \frac{\partial^2 u}{\partial x^2}, \dots, p_1(x), p_2(x), \dots U(t)\right),
\end{align}
with the initial condition $u(x,0)=u_0(x)$, where $ F $ is the PDE that models the system dynamics, $ U(t) $ is the control function, and $ p_i(x) $, for $ i = 1, \dots, m $ are the coefficient functions. The PDE control problem is to design a controller which stabilizes the system or regulates states to a desired trajectory.

\subsection{Backstepping control operator}
Employing the PDE backstepping method \cite{krstic2008boundary}, a feedback controller—typically implemented as a boundary controller—can be designed in the following form:
\begin{align} 
	\label{U(t)}
	U(t) &= \int_{0}^{1} k(1,y)u(y, t) \, dy,
\end{align}
where $k(1,y)$ is the control gain derived from the backstepping kernel function. This kernel depends on the system coefficient functions $p_i(x)$, $i=1,\ldots, m$.
To formalize the controller structure, we define the control operator $\mathcal{U}: (C[0,1])^m\times L^2[0,1] \mapsto \mathbb{R}$, which maps the coefficient functions $(p_1, \ldots, p_m) \in (C[0,1])^m $ and the system state $u \in L^2[0,1] $ to the control input $U $. Specifically,
\begin{equation}
    U := \mathcal{U}(p_1, \dots, p_m, u),
\end{equation}
where $U $ is defined as in \eqref{U(t)}. 

\subsection{DeepONet for learning  backstepping controller}\label{subsec_DeepONet}
DeepONet is a deep learning architecture designed to learn function-to-function mappings, consisting of branch and trunk networks. The output function is a weighted sum of basis functions, with the branch network learning the weights and the trunk network learning the basis functions.



For the PDE controller approximated by a DeepONet, the DeepONet receives coefficient functions $p_i(x)$ and state $u(x,\cdot)$ as inputs, generating a scalar control signal $U(\cdot)$ as output. In the mapping, the input and output are time-independent, the control input is generated based on the current system state. All training data are generated by solving the PDE \eqref{PDE} with the backstepping controller \eqref{U(t)}, using different initial conditions and coefficient functions randomly sampled from Chebyshev polynomials with random parameters \cite{bhan2023neural}.

To train the DeepONet, we employ the following state-based loss function:
\begin{align} 
	\label{Loss function}
\varepsilon = |U_{NO}-U|,
\end{align}
where \(U_{NO}\) denotes the control generated by DeepONet, while \(U\) denotes the backstepping control.

\begin{figure}[htbp]
	\centering
	\includegraphics[width=0.8\linewidth]{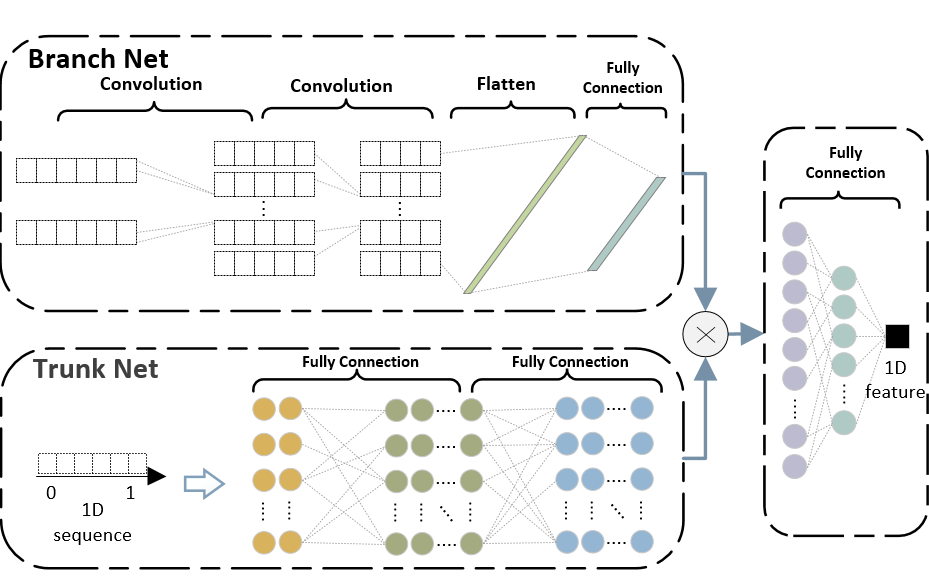}
	\caption{\small Architecture of DeepONet as feature extractor for SAC.}
	\label{deeponet}
\end{figure}

To integrate the DeepONet into the reinforcement learning architecture, we modify its output from a scalar control signal to a vector representation, enabling feature extraction for both actor-critic networks in SAC, as depicted in Fig. \ref{deeponet}. The dimensionality of the DeepONet output is determined by its branch and trunk networks, along with subsequent fully connected layers. The pre-trained DeepONet encapsulates the knowledge of the backstepping controller, enabling seamless integration of prior control strategies into the RL decision-making process. This incorporation not only improves the training efficiency of the reinforcement learning algorithm but also enhances the performance of the RL-based controller. 

\section{NO-based reinforcement learning framework}
\begin{figure}[htbp]
	\centering
	\includegraphics[width=0.8\linewidth]{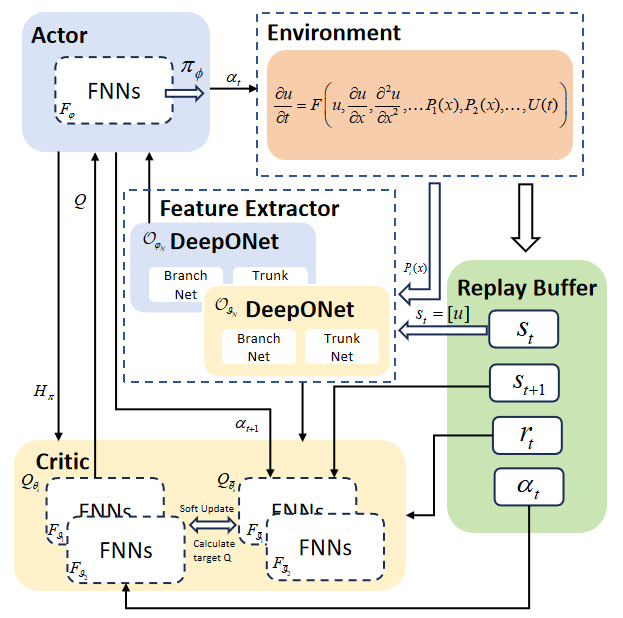}
	\caption{\small The DeepONet pre-trained with the backstepping method is embedded into the SAC framework.}
	\label{fig:1}
\end{figure}

In this section, we elaborate on the reinforcement learning (RL) framework that fuses a pre-trained DeepONet with the Soft Actor-Critic (SAC) algorithm. This framework is composed of an actor network $\pi_\phi(a_t|s_t) $, a pair of critic networks denoted as $Q_{\theta_i}(s_t,a_t) $, as well as their respective target networks $Q_{\bar\theta_i}(s_t,a_t) $. Serving as a feature extraction component, the DeepONet replaces the CNNs adopted in the canonical SAC structure, and establishes a direct connection with the fully connected layers of the actor and critic networks. In the training phase, the actor-critic networks and the embedded DeepONet feature extractor are optimized in a joint manner through the backpropagation algorithm. The schematic diagram of this RL framework is presented in Fig. \ref{fig:1}.

\subsection{Actor network} By introducing the policy entropy $ H_{\pi} $, it quantifies the uncertainty of the action distribution conditioned on state $ s_t $. The policy entropy is defined as: $H(\pi(\cdot|s_t)) = -\mathbb{E}_{a_t \sim \pi_{\phi}} [\log \pi_{\phi}(a_t|s_t)]$. The policy network $\pi_{\phi}(a_t|s_t)$ is trained by optimizing the following subsequent objective function:
\begin{align} 
	\label{J}
	J_\pi({\phi}) = \mathbb{E}_{s_t \sim \mathcal{D}} \mathbb{E}_{a_t \sim \pi_{\phi}} \left[ -\alpha H(\pi_{\phi}(\cdot|s_t)) - \min_{i=1,2} Q_{\theta_i}(s_t, a_t) \right],
\end{align}
with $\mathcal{D}$ characterizes the distribution of system states. Consequently, the update of the policy network parameter $\phi$ is realized by
\begin{align} 
	\label{pi}
\phi \leftarrow \phi - \eta \nabla_\phi J_\pi(\phi).
\end{align}

\subsection{Critic networks} We use a double Q-network structure (two independent FNNs), the action-value function is defined as $ Q_\pi(s_t, a_t) $, which represents the expected cumulative return that the agent obtains by following the policy $ \pi_\phi $ after executing the action $ a_t $ in the state $ s_t $. Using two target networks $Q_{\bar\theta_{1,2}}$ to determine the target Q-function:

\begin{align} 
	\label{yt}
	y_t = r_t + \gamma \left( \min_{i=1,2} Q_{\bar\theta_i}(s_{t+1}, \hat{a}_{t+1}) - \alpha \log \pi_{\phi}(\hat{a}_{t+1}|s_{t+1}) \right),
\end{align}
$\alpha$ is the temperature parameter for balancing reward maximization and policy entropy, and the action $\hat{a}_{t+1} \sim \pi(\cdot|s_{t+1})$ is drawn from the actor network. The weights of the target network $Q_{\bar\theta_i} $ are updated by
\begin{align} 
\label{-theta}
\bar{\theta}_i \leftarrow \tau\theta_i + (1-\tau)\bar{\theta}_i,
\end{align}
where $ i \in \{1, 2\} $ and $ \tau \in (0, 1] $ denotes the weighted coefficient.

The action-value networks $ Q_{\theta_1}(s_t, a_t) $ and $ Q_{\theta_2}(s_t, a_t) $ are trained by minimizing the
following objective:
\begin{align} 
	\label{JQ}
	J_Q(\theta_i) = \mathbb{E}_{(s_t,a_t) \sim \mathcal{W}} \left[ (Q_{\theta_i}(s_t, a_t) - y_t)^2 \right],
\end{align}
where $ i \in \{1, 2\} $ and $\mathcal{W}$ is the distribution of the states and actions. The parameter $\theta_i$ is updated by stochastic gradient descent:
\begin{align}
\label{Jtheta}
\theta_i \leftarrow \theta_i - \lambda \hat\nabla_{\theta_i} {J}_Q(\theta_i).
\end{align}

\subsection{Algorithm (Markov Decision Process)}
We have briefly outlined the processing of the Markov Decision Process (MDP), as shown in Algorithm \ref{alg:deepsac}.

1) State space $ \mathcal{S} $: we denote $ s_t = \{u(x,t)\} \in \mathcal{S} $, which constitutes the extended state at time $ t $.

2) Action space $ \mathcal{A} \subseteq \mathbb{R} $: it is $ \mathcal{A} = [-\bar{U},\bar{U}] $, with $ \bar{U} = \sup_{t \in \mathbb{R}^+} |U(t)| $, $U(t) = a_t \in \mathcal{A}$.

3) Policy function $\pi(a_t|s_t)$: this policy is formulated as a probability density function $\pi(a_t|s_t)$, which establishes a mapping relationship from the current state space $\mathcal{S}$ to the probability distribution over the action space $\mathcal{A}$.

4) State transition $ p(s_{t+1} | s_t, a_t) $: this $ p(s_{t+1} | s_t, a_t) $ is the probability distribution of $ s_{t+1} $, which depends on the current state $ s_t $ and the action $ a_t $ chosen by the agent.

5) Reward function $r_t$ and additional reward function $q$: it is designed to guide the learning of the optimal control policy. The reward function is defined as follows:
\begin{align}
r_t(s_t, s_{t+1}) &= -1 \cdot \| s_{t+1} - s_t \|_{L_2}, \label{eq:reward} \\
q(s_T, a_0, \ldots, a_T) &= 
\begin{cases} 
0 &  \| s_T \|_{L_2} > \zeta, \\
\begin{aligned}
\sigma - \dfrac{1}{\eta} \cdot \sum_{\tau=0}^T |a_\tau|_{L_1} \\
- \| s_t \|_{L_2} 
\end{aligned} &  \| s_T \|_{L_2} \leq \zeta,
\end{cases} \label{eq:q_function}
\end{align}
here, $\sigma$, $\eta$, and $\zeta$ are hyperparameters. The reward function $r_t $ drives state convergence at each step, with an additional reward $q $ given at the end of each episode if the $L_2 $ norm of the final state is less than or equal to the threshold $\zeta $.

 \begin{algorithm}
	\caption{DeepONet-SAC Training Process}
	\label{alg:deepsac}
	\begin{algorithmic}[1]
		\State Initialize network parameters: $\phi, \theta_i, \bar{\theta}_i$ ($i=1,2$), replay memory $\mathcal{W}$, and starting state $s_0 = \{u_0(x)\}$
		\For{iteration $l$ from $0$ to $\text{max steps}$}
		\If{$\|s_t\|_{L_2} > \text{threshold}$ \textbf{or} $t \geq T$}
		\State Reinitialize state $s_t \leftarrow s_0$
		\EndIf
		\State Set the instantaneous state $s_t = \{u_t\}$
		\State Draw action $a_t \sim \pi_\phi(\cdot \mid s_t)$
		\State Apply $a_t$, observe next state $s_{t+1}$ and reward $r_t$
		\State Store transition $(s_t, a_t, r_t, s_{t+1})$ in $\mathcal{W}$
		\For{each gradient update iteration}
		\State Randomly sample mini-batch from $\mathcal{W}$
		\State Calculate target value $y_t$ via equation (6)
		\State Adjust critic parameters $Q_{\theta_{1,2}}$ using (8)
		\State Update target network parameters $Q_{\bar{\theta}_{1,2}}$ by (5)
		\State Modify actor parameters $\phi$ according to (4)
		\EndFor
		\EndFor
	\end{algorithmic}
\end{algorithm}

\section{Simulation} 
In this section, the proposed RL controller is applied to both hyperbolic and parabolic PDEs, and the details of the backstepping design can be found in \cite{krstic2008boundary}. We adopt the standard SAC implementation from the PDE Control Gym platform introduced in \cite{bhan2024pde}. Our proposed method—SAC embedded with DeepONet pre-trained using the backstepping controller (referred to as NOSAC\_training)—is compared against four baselines: the backstepping controller (Backstepping), standard SAC (SAC), and SAC embedded with DeepONet without pretraining (NOSAC).

Experiments conducted on an Intel i9-13900KF/RTX 4090 workstation; SAC and reward hyperparameters in \cite{bhan2024pde} Appendix A.2.

\subsection{Simulation of a 1D Hyperbolic PDE}
First, A DeepONet is trained to emulate the backstepping controller, with corresponding training data generated accordingly.
We employ the finite difference method with $dx = 0.01$ over the time interval $[0, 5]$ with a temporal step size of $\Delta t = 0.001$. The coefficient $\beta(x)$ is sampled from the Chebyshev polynomial $\beta(x) = 5 \cos\big(\gamma \cos^{-1}(x)\big)$, where $\gamma \sim U[5.5, 7]$. The initial condition $u_0(x)$ is assumed to be a constant randomly chosen from the uniform distribution $U[1, 10]$. 

We generate a dataset of $6 \times 10^5$ samples, which comprises 100 different $\beta(x)$ and 60 different $u_0(x)$. We collect data every 50 time steps, and then randomly shuffle and split the dataset into a training set 90\% and a testing set 10\%. The training process requires approximately 175 minutes.

\begin{figure}[htbp]
	\centering
	\subfloat[]{\epsfig{figure=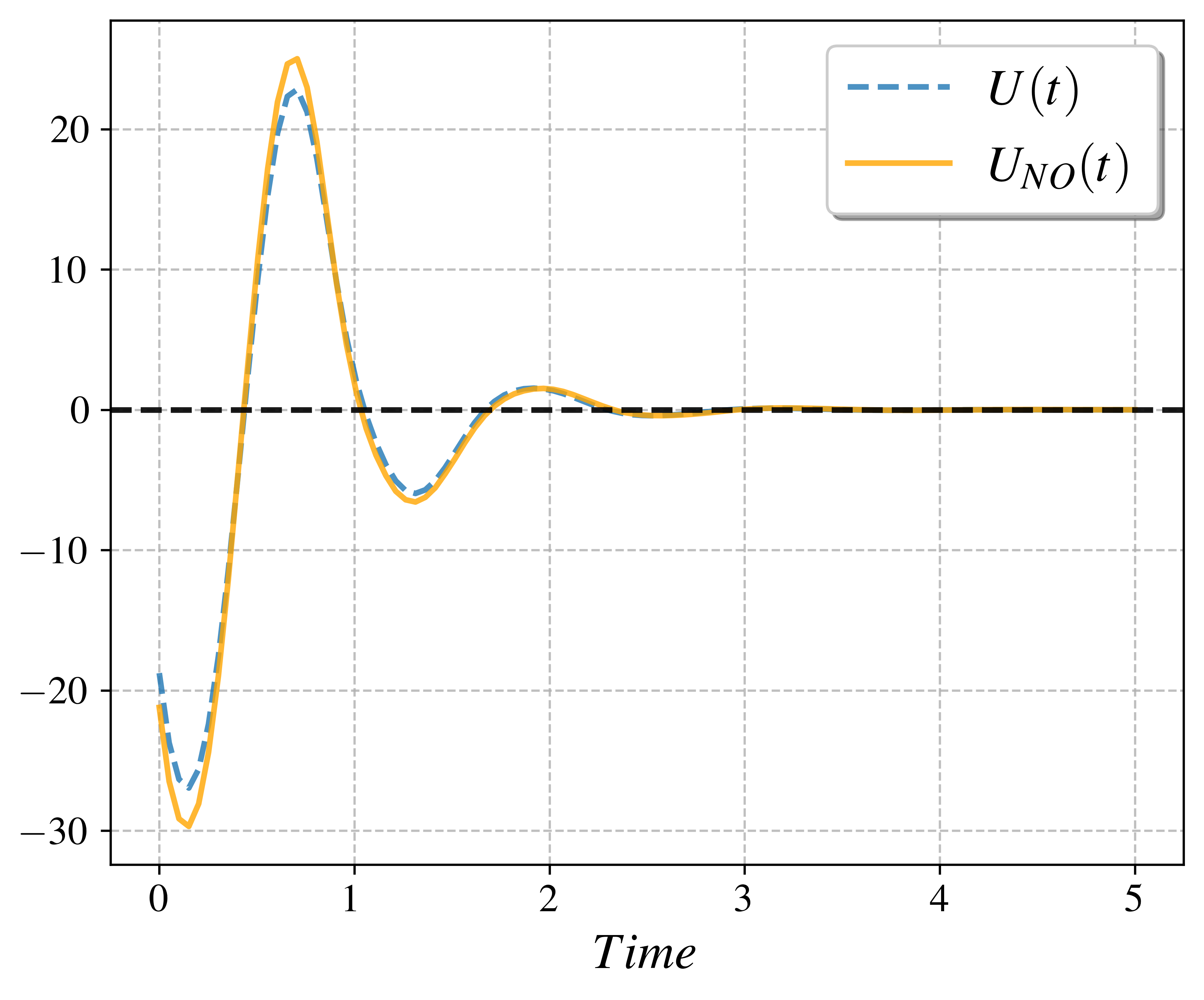, width=0.2\textwidth}}
	\subfloat[]{\epsfig{figure=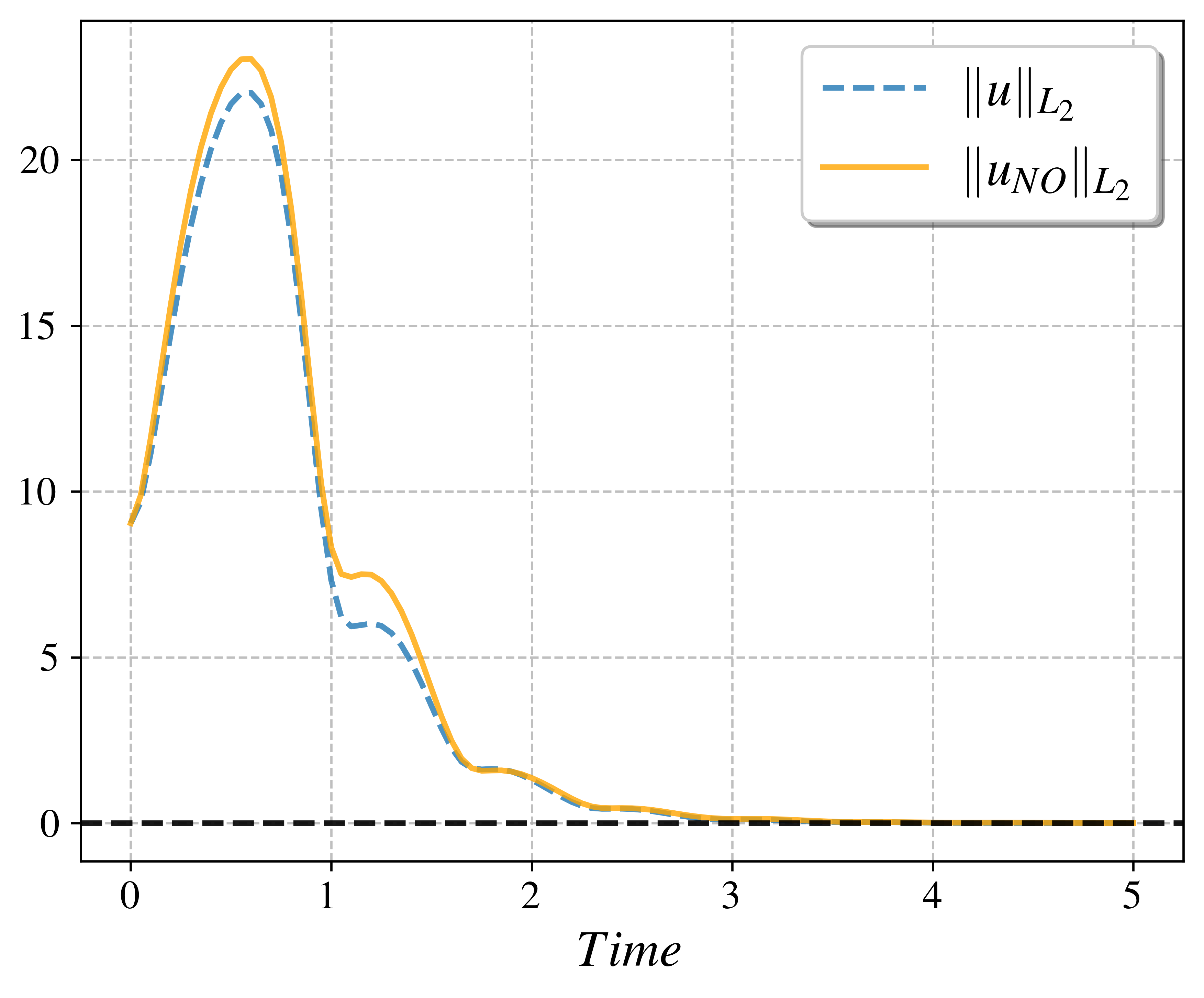, width=0.2\textwidth}}

	\caption{\small (a) Control input and (b) $L_2$ norm of the state $u(x,t)$ for hyperbolic PDE with $\gamma = 5.5$ and $u_0(x) = 9$ under the backstepping controller and the DeepONet-approximated controller.}
	\label{fig:3}
\end{figure}

The DeepONet approximation results are shown in Fig. \ref{fig:3}, illustrating the boundary control $U(t)$ and the $L_2$ norm of the state under the backstepping controller and the NO-based controller for $\gamma = 5.5$. It can be seen that $U_{NO}$ remains almost indistinguishable from $U$, and the state converges to zero under both controllers ($\|u_{NO} - u\|_{L_2}$ is less than $10^{-4}$ after about $3$ seconds).
\begin{figure}[htbp]
	\centering
	\subfloat[]{\epsfig{figure=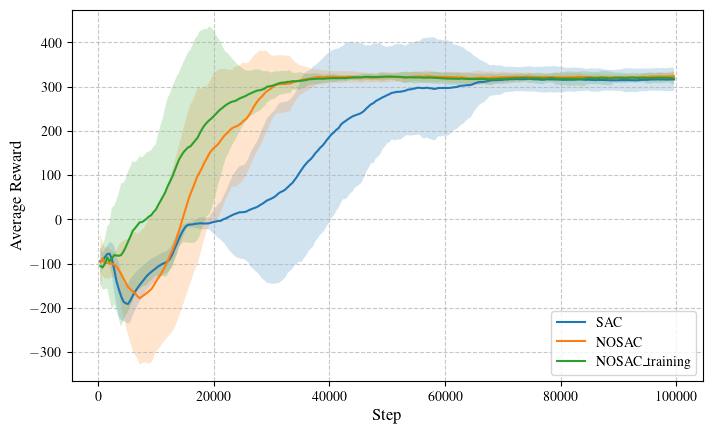, width=0.24\textwidth}}
	\subfloat[]{\epsfig{figure=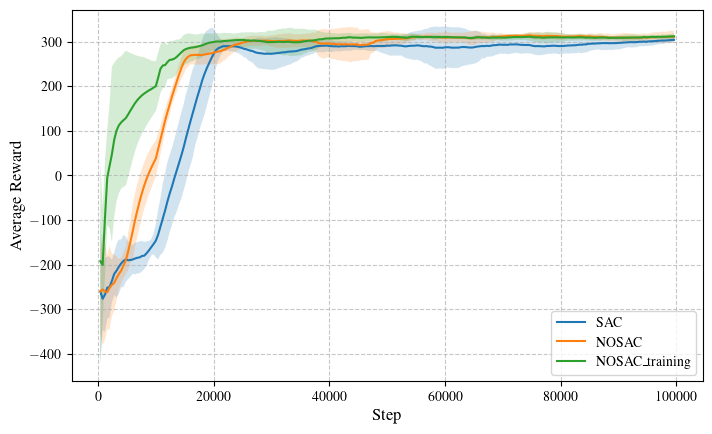, width=0.24\textwidth}}
	
	\caption{\small Reward curves for the three RL training processes: (a)  1D Hyperbolic PDE, (b)  1D Parabolic PDE.}
	\label{jl}
\end{figure}

In the RL training, we set $\gamma=5.5$ and sample $u_0(x)$ sampled from a uniform distribution $U[1, 10]$ at each iteration. The training consists of 100{,}000 steps with 100 interaction steps per episode. The training time is about 18, 20, and 11 minutes for NOSAC\_training, NOSAC, and SAC taking, respectively. Fig. \ref{jl}(a) shows NOSAC\_training achieves faster reward growth and policy convergence compared to the others, thanks to the pre-trained DeepONet. Pre-trained on backstepping controller, the DeepONet captures stable control behaviors, enabling it to encode this knowledge into feature vectors with higher fidelity. This serves as a warm start, for the SAC actor—reducing exploration costs and increasing initial rewards—and provides the critic with more informative state features.

\begin{figure}[htbp]
	\centering  
	
	\subfloat[]{\includegraphics[width=0.24\textwidth]{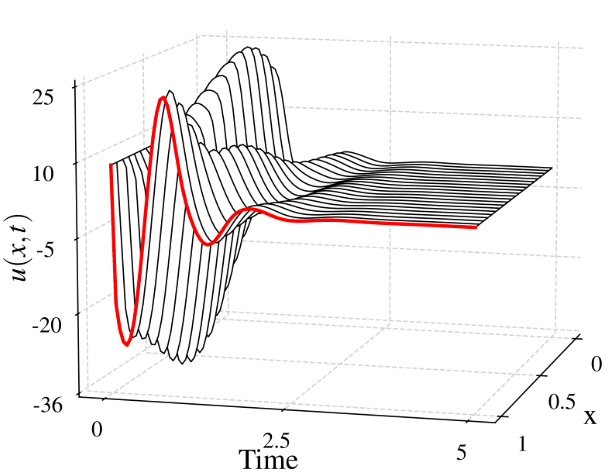}}
	\hfill  
	\subfloat[]{\includegraphics[width=0.24\textwidth]{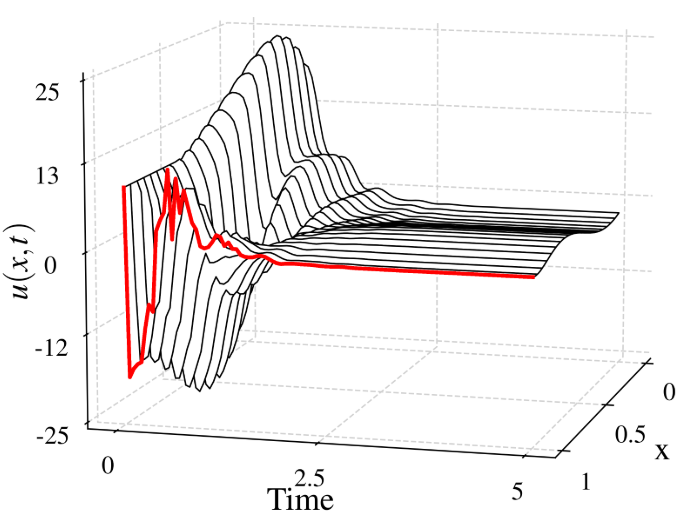}}

	\subfloat[]{\includegraphics[width=0.24\textwidth]{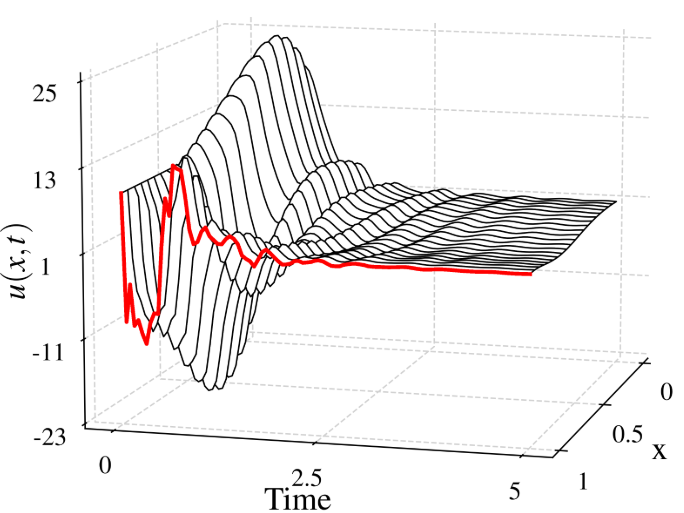}}
	\hfill  
	\subfloat[]{\includegraphics[width=0.24\textwidth]{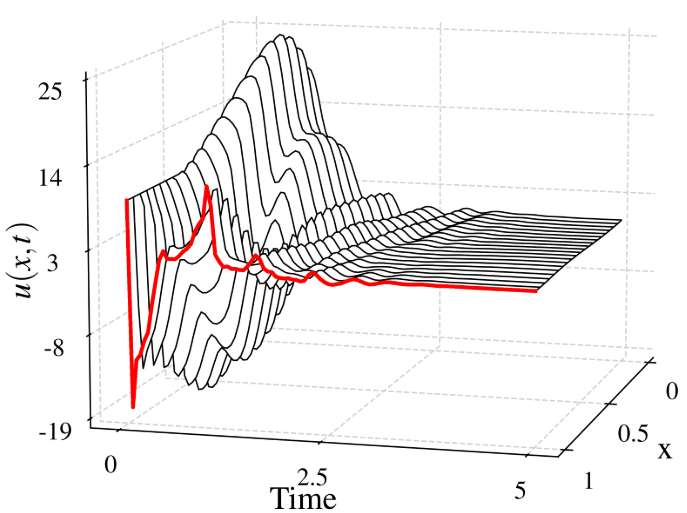}}
	\caption{\small 
		Closed-loop state evolution with $\gamma = 5.5$ and $u_0(x) = 9$ under         (a) Backstepping,
		(b) SAC,
		(c) NOSAC, 
		(d) NOSAC\_training.
	}
	\label{fig:5}
\end{figure}

\begin{figure}[htbp]
	\centering
	\subfloat[]{\epsfig{figure=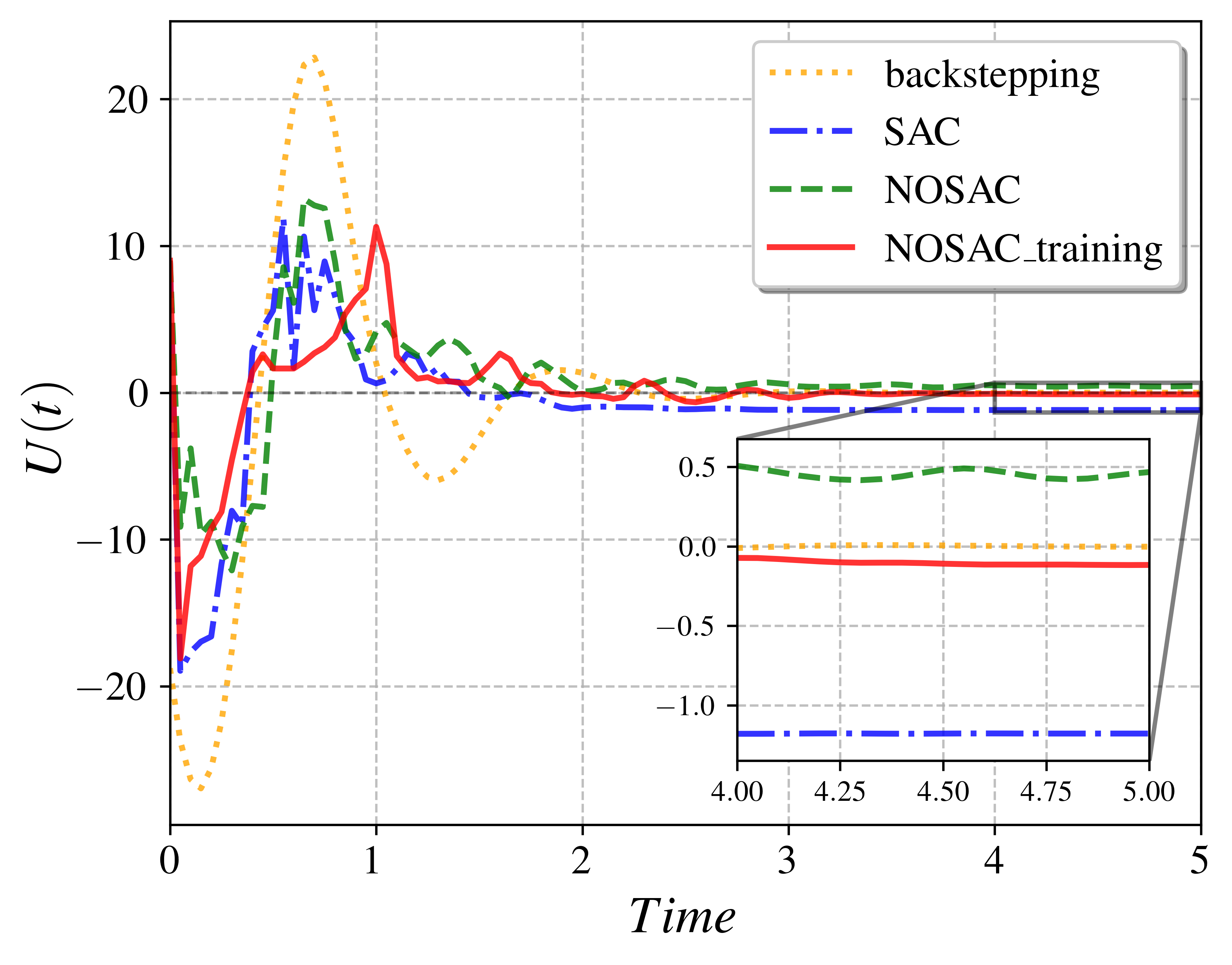, width=0.24\textwidth}}
	\subfloat[]{\epsfig{figure=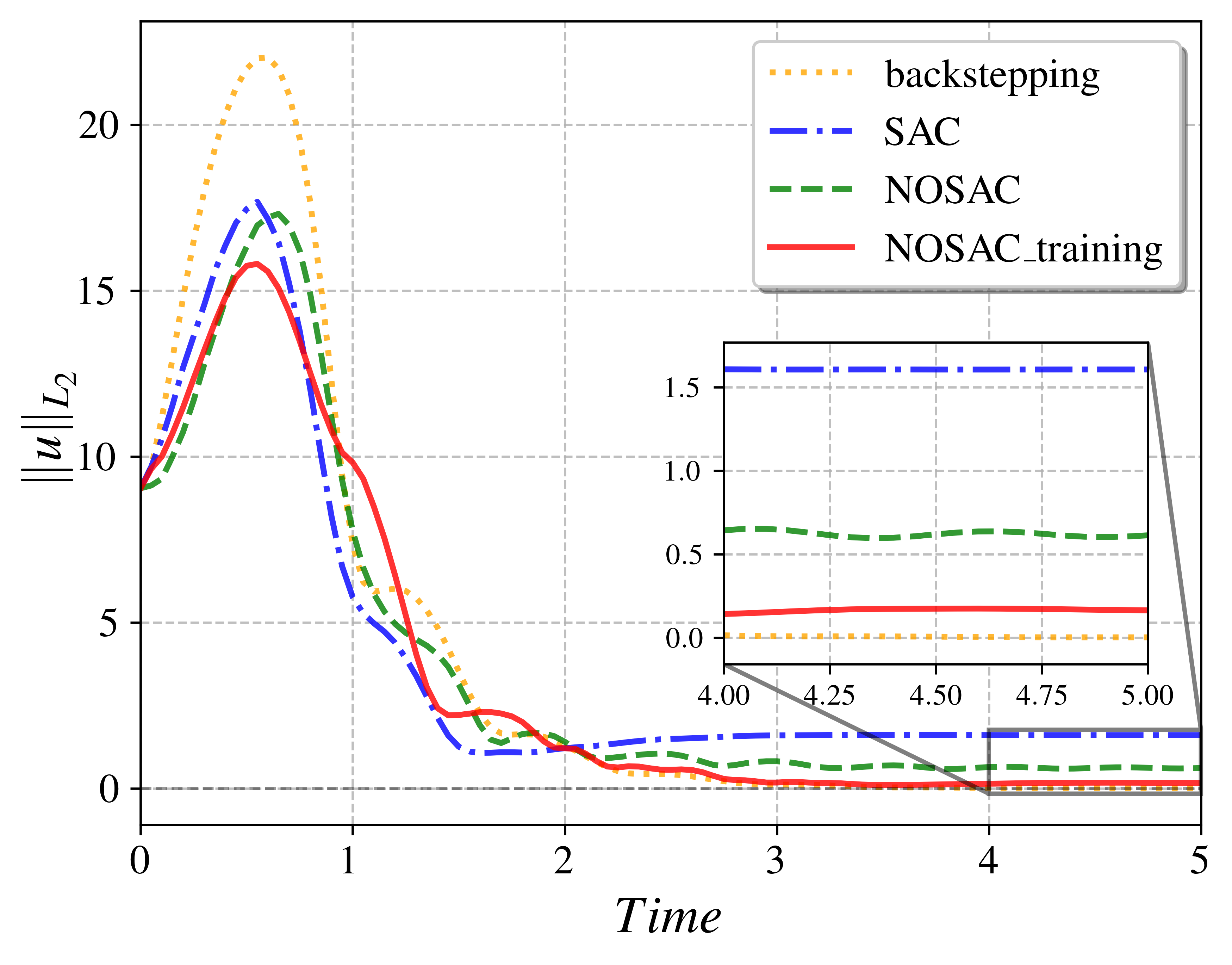, width=0.24\textwidth}}
	
	\caption{\small (a) Control input and (b) $L_2$ norm of the state with $\gamma = 5.5$ and $u_0(x) = 9$ under the four different controllers. }
	\label{fig:6}
\end{figure}

We apply the NOSAC\_training controller to the 1D hyperbolic PDE and compare its performance with the backstepping, SAC, and NOSAC, as shown in Fig. \ref{fig:5}. The backstepping controller yields smoother closed-loop state evolution, while the proposed NOSAC\_training controller exhibits less overshoot. To better illustrate performance, Fig. \ref{fig:6} plots the control effort and the $L_2$ norm of the state, showing that NOSAC\_training outperforms the other controllers in terms of overshoot and convergence rate, whereas the backstepping controller demonstrates superior steady-state precision.

\begin{figure}[htbp]
	\centering
	\subfloat[]{\epsfig{figure=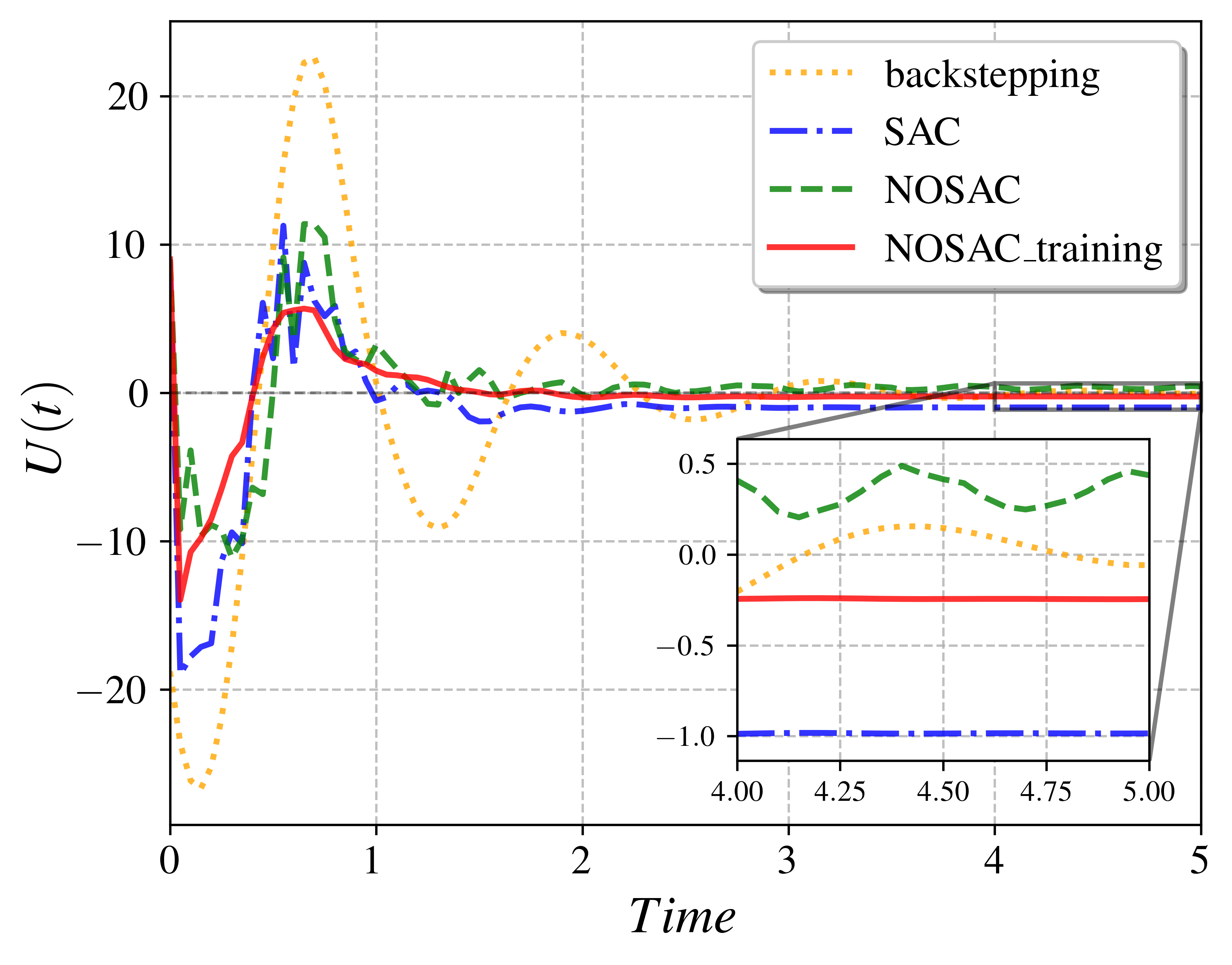, width=0.24\textwidth}}
	\subfloat[]{\epsfig{figure=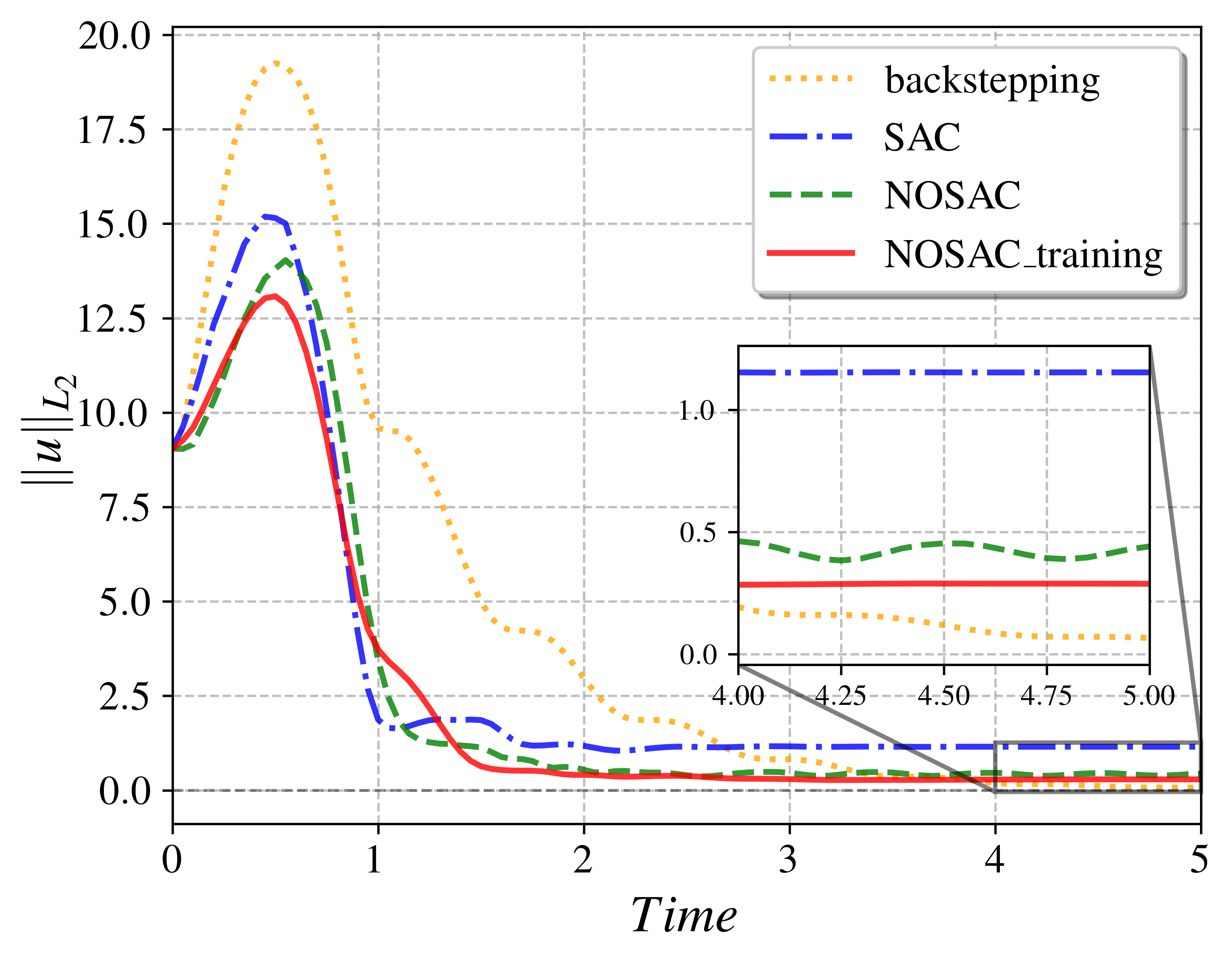, width=0.24\textwidth}}
	
	\caption{\small 
		 (a) Control input and (b) $L_2$ norm of the state under model mismatch with system's parameter $\gamma = 5.7$ and initial condition $u_0(x) = 9$ while training the three RLs and designing the backstepping using  $\gamma = 5.5$. }
	\label{fig:8}
\end{figure}

To evaluate robustness under model mismatch, we simulate a system with $\gamma = 5.7$ using the three LR\_based controllers trained at $\gamma = 5.5$ and the backstepping controller designed for $\gamma = 5.5$. The results in Fig. \ref{fig:8} demonstrate that NOSAC\_training achieves superior robustness, outperforming the SAC and NOSAC controllers in overshoot, convergence speed, and steady-state error.


\subsection{Simulation of a 1D Parabolic PDE}
When dealing with 1D Parabolic PDE, the implicit method is employed, which enhances accuracy and reduces computation. We employ the $dx = 0.01$ over the time interval $[0,1]$ with a temporal step size of $\Delta t = 0.001$. The coefficient $\lambda(x)$ is sampled from the Chebyshev polynomial $\lambda(x) = 50 \cos\big(\gamma \cos^{-1}(x)\big)$, where $\gamma \sim U[8, 12]$. We generate a dataset comprising $6 \times 10^5$ samples, which were produced by combining 100 different $\lambda(x)$ and 60 different $u_0(x)$. During the data collection process, data was collected every 100 time steps, its DeepONet training process is the same as that for the 1D Hyperbolic PDE. The total training process takes about 190 minutes.

\begin{figure}[htbp]
	\centering
	\subfloat[]{\epsfig{figure=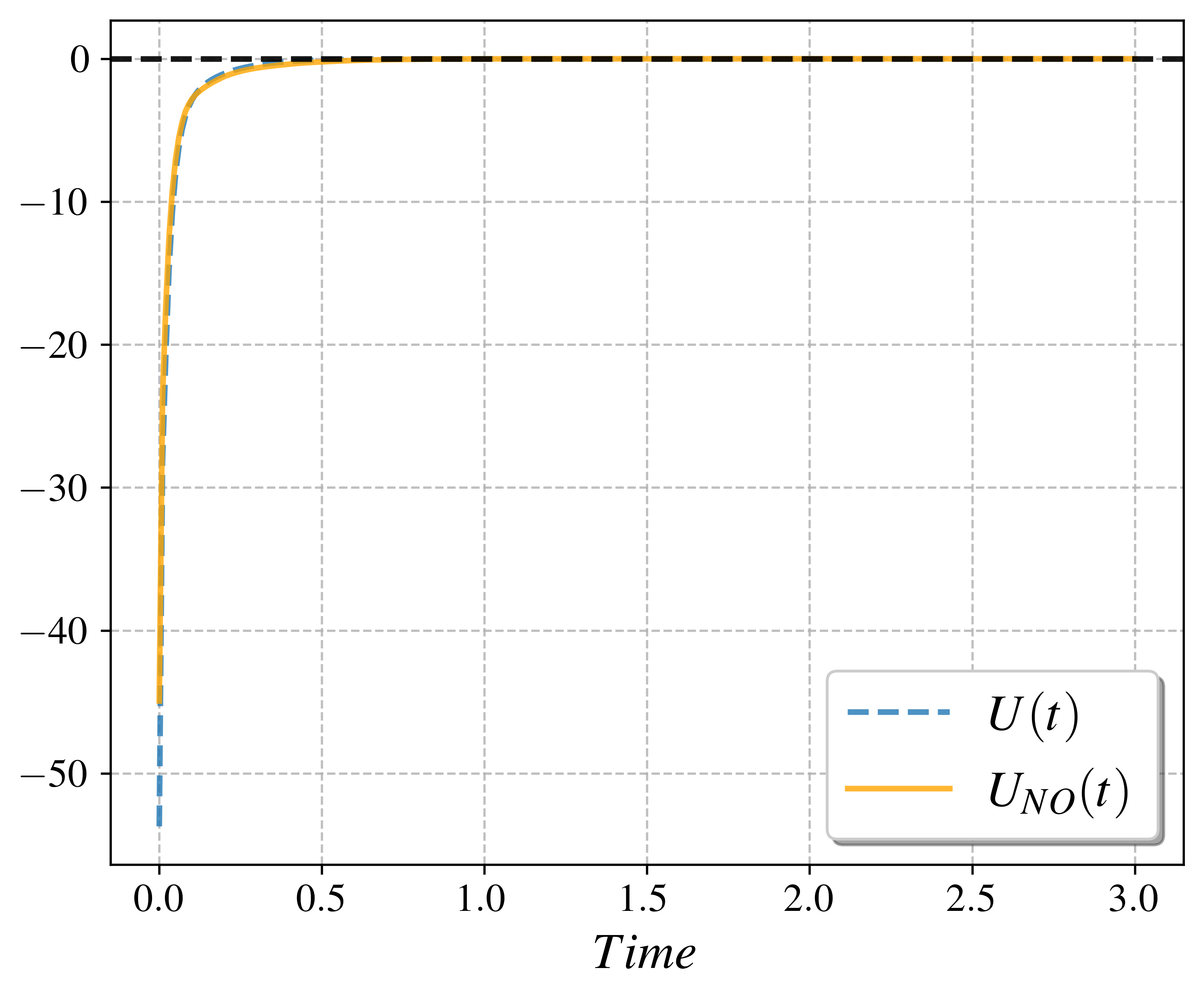, width=0.2\textwidth}}
	\subfloat[]{\epsfig{figure=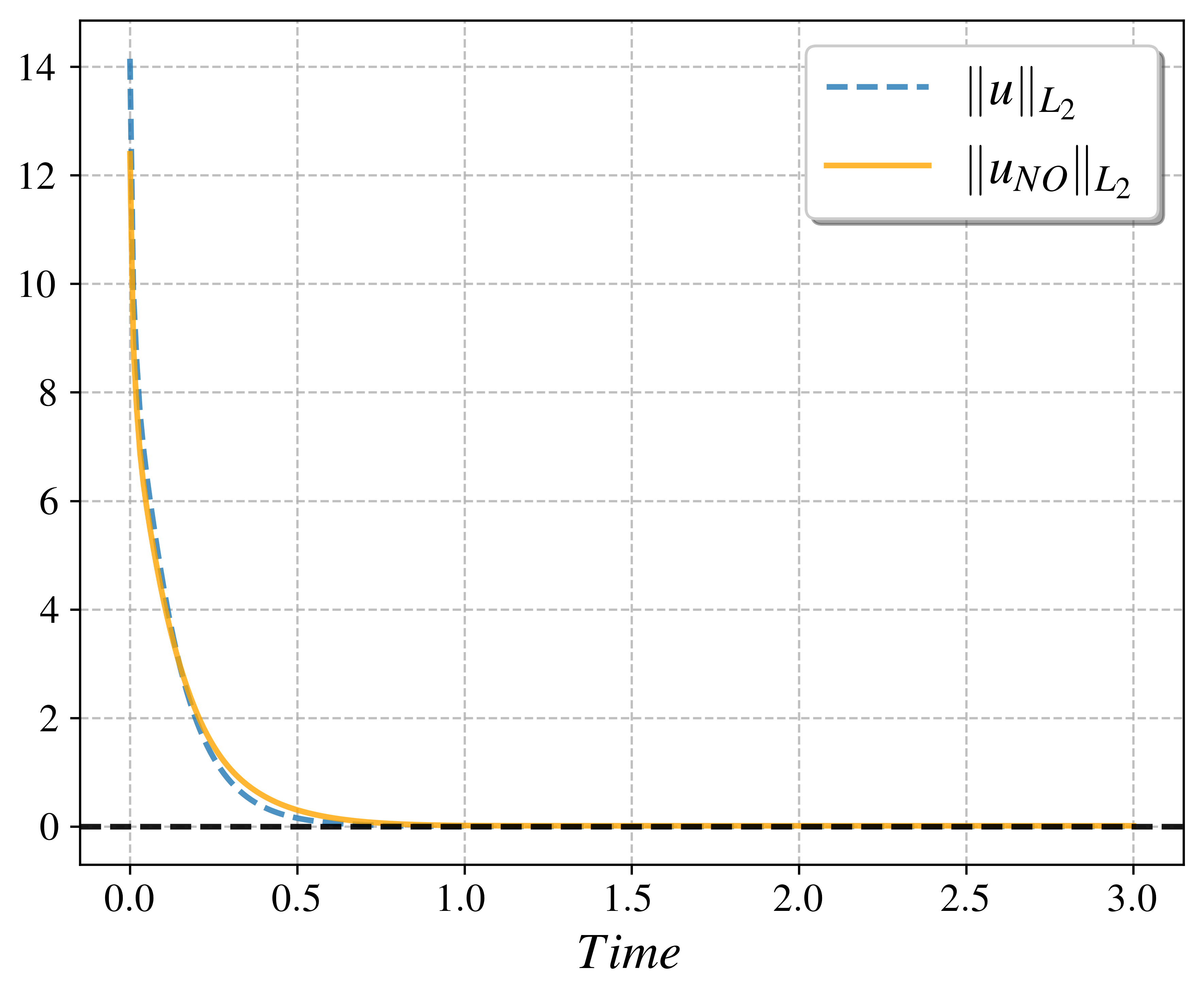, width=0.2\textwidth}}
	
	\caption{\small (a) Control input and (b) $L_2$ norm of the state $u(x,t)$ for parabolic PDE with $\gamma = 9$ and $u_0(x) = 9$ under the backstepping controller and the DeepONet-approximated controller.}
	\label{fig:10}
\end{figure}

The DeepONet approximation results are shown in Fig. \ref{fig:10}, illustrating the boundary control and the $L_2$ norm of the state under the backstepping controller and the NO-based controller for $\gamma = 9$. It can be seen that $U_{NO}$ remains almost indistinguishable from $U$, and the state converges to zero under both controllers($\|u_{NO} - u\|_{L_2}$ is less than $10^{-4}$ after about $0.5$ seconds). 

In the RL training, we set $\gamma=9$, while the other settings are consistent with those used in the hyperbolic PDE. The training processes of NOSAC\_training, NOSAC, and SAC take approximately 20 minutes, 22 minutes, and 14 minutes, respectively. Fig. \ref{jl}(b) shows the reward curves for the proposed method and other two baseline algorithms on 1D Parabolic PDE, demonstrating that NOSAC\_training exhibits the fastest reward growth and convergence.

\begin{figure}[htbp]
	\centering  
	
	\subfloat[]{\includegraphics[width=0.24\textwidth]{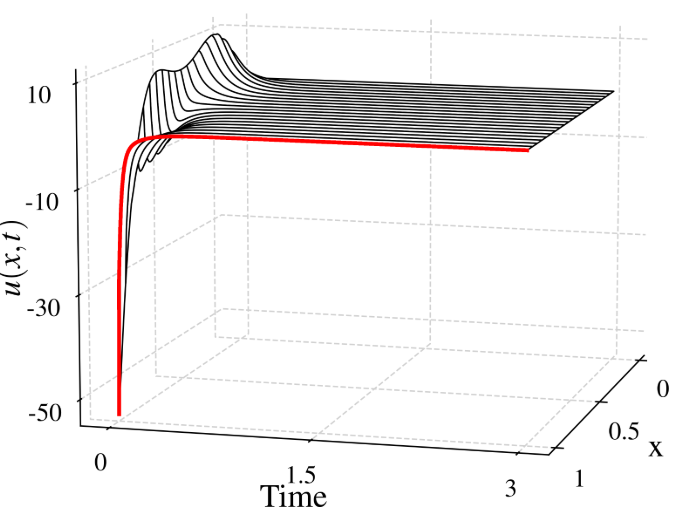}}
	\hfill  
	\subfloat[]{\includegraphics[width=0.24\textwidth]{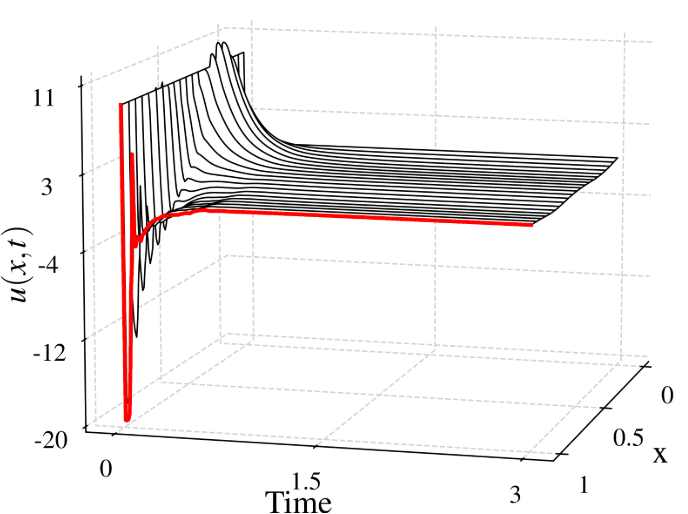}}

	\subfloat[]{\includegraphics[width=0.24\textwidth]{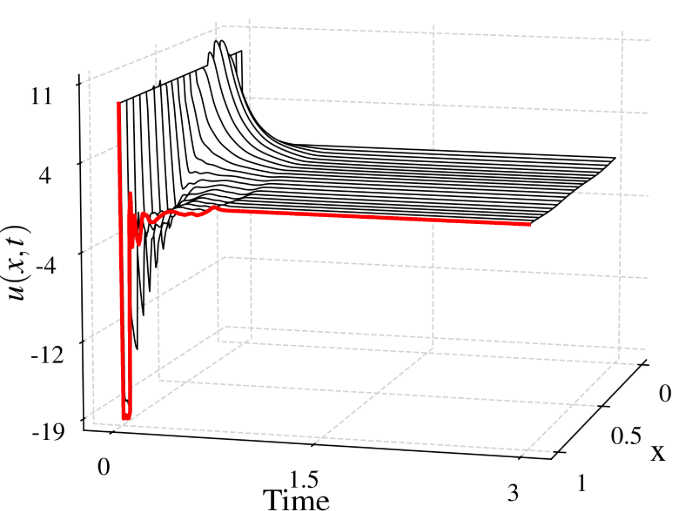}}
	\hfill  
	\subfloat[]{\includegraphics[width=0.24\textwidth]{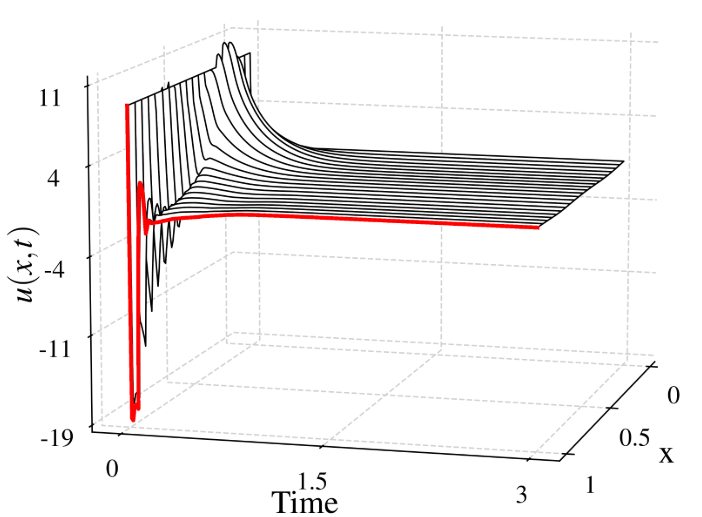}}
	\caption{\small 
		\small 
		Closed-loop state evolution with $\gamma = 9$ and $u_0(x) = 9$ under         (a) Backstepping,
		(b) SAC,
		(c) NOSAC, 
		(d) NOSAC\_training.
	}
	\label{fig:z}
\end{figure}

\begin{figure}[htbp]
	\centering
	\subfloat[]{\epsfig{figure=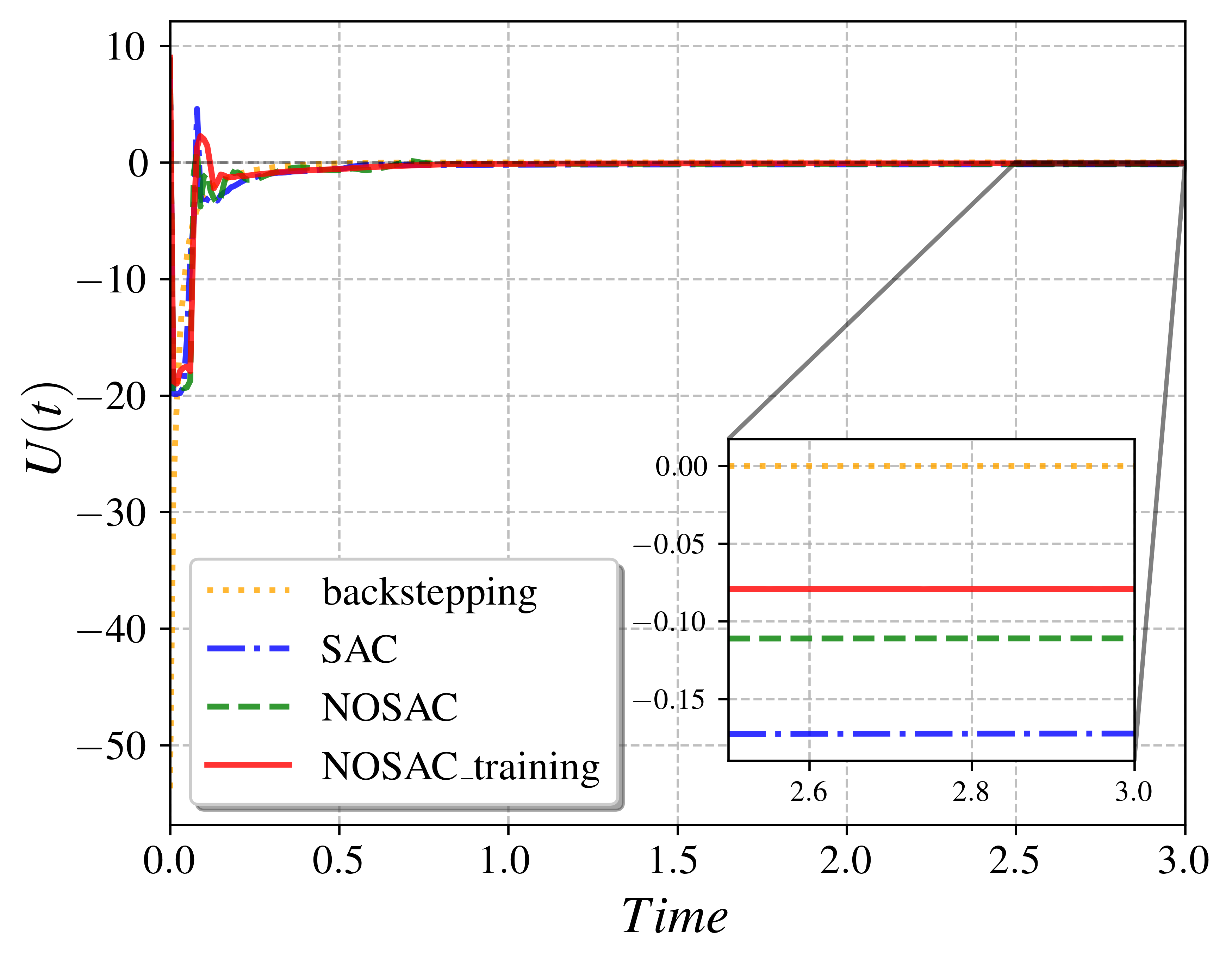, width=0.24\textwidth}}
	\subfloat[]{\epsfig{figure=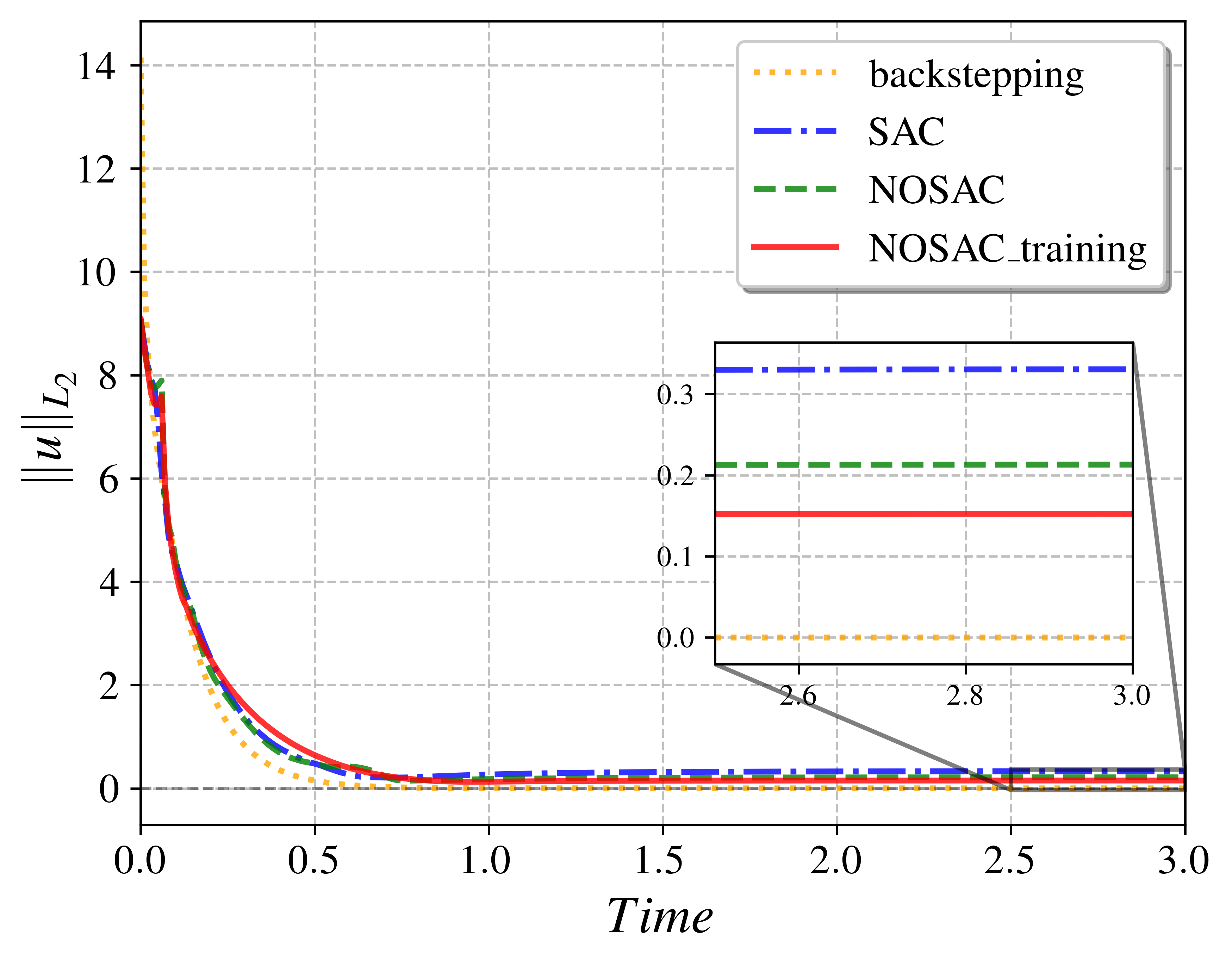, width=0.24\textwidth}}
	
	\caption{\small (a) Control input and (b) $L_2$ norm of the state with $\gamma = 9$ and $u_0(x) = 9$ under the four different controllers.}
	\label{fig:13}
\end{figure}

We apply the NOSAC\_training controller to the 1D parabolic PDE and compare its performance with the backstepping, SAC, and NOSAC, as shown in Fig. \ref{fig:z}. The proposed NOSAC\_training controller exhibits less overshoot. Fig. \ref{fig:13} plots the control effort and the $L_2$ norm of the state, showing that NOSAC\_training is optimal in terms of overshoot and convergence speed, conversely, the backstepping controller exhibits enhanced steady-state convergence.

\begin{figure}[htbp]
	\centering
	\subfloat[]{\epsfig{figure=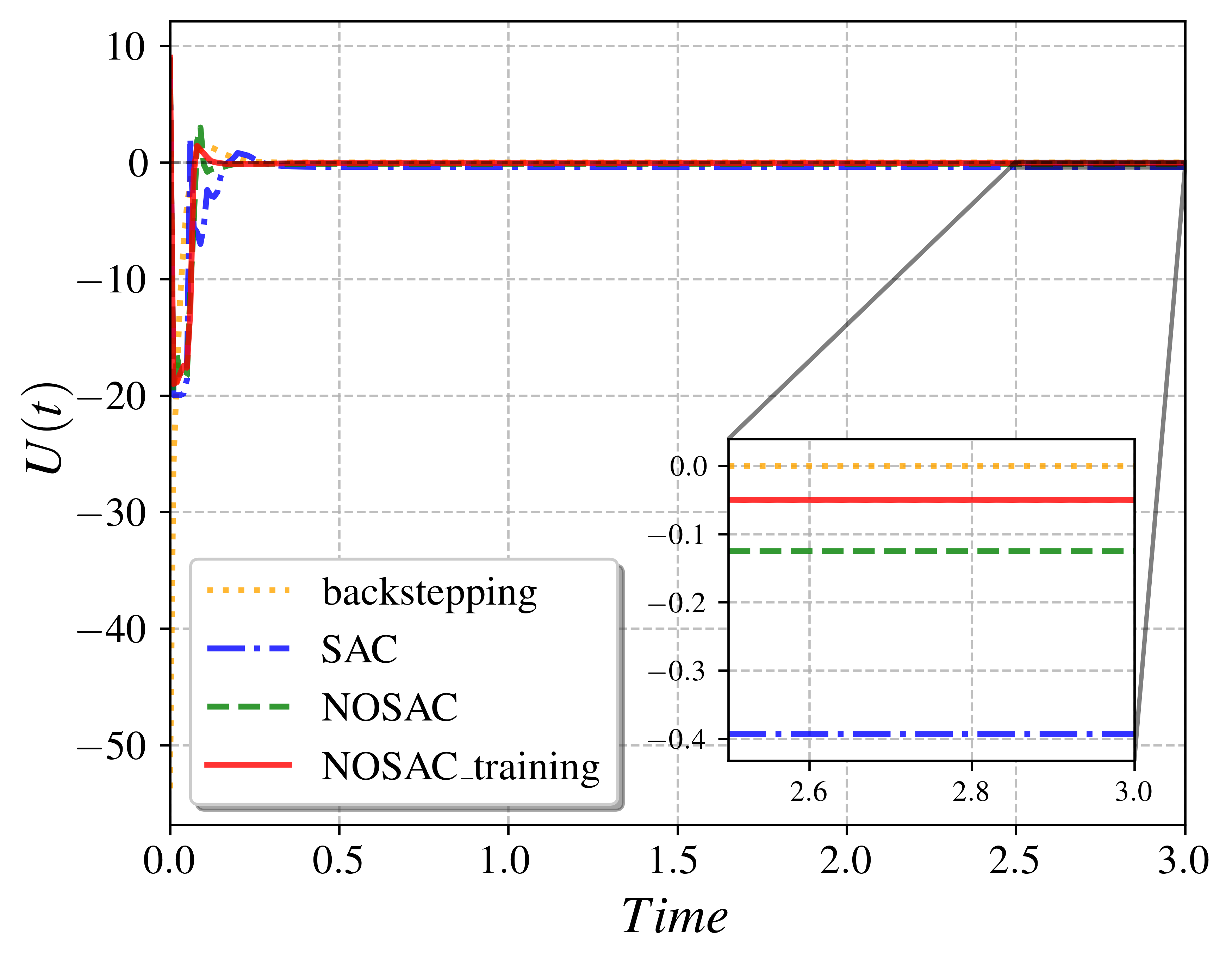, width=0.24\textwidth}}
	\subfloat[]{\epsfig{figure=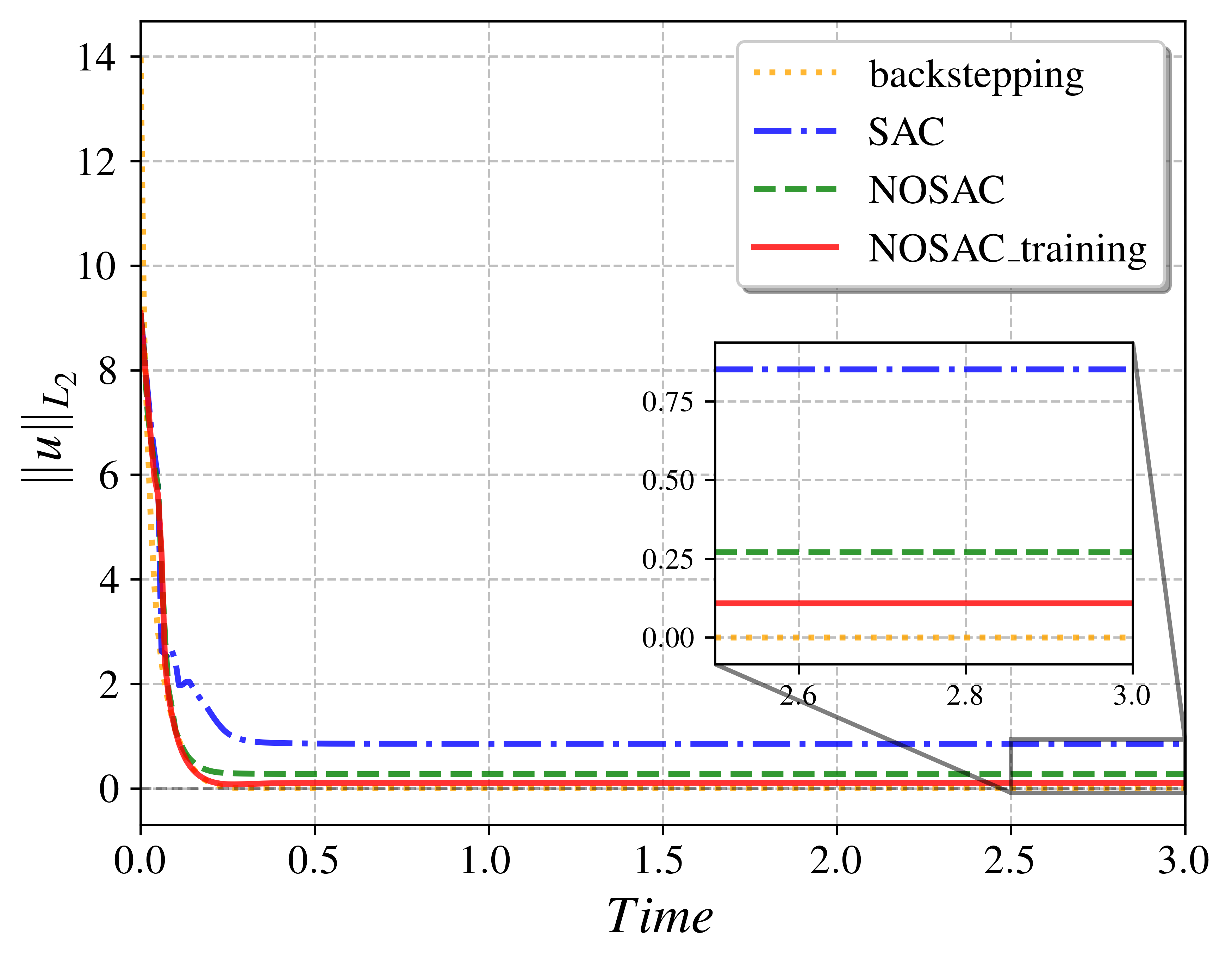, width=0.24\textwidth}}
	
	\caption{\small 
		 (a) Control input and (b) $L_2$ norm of the state under model mismatch with system's parameter $\gamma = 8.5$ and initial condition $u_0(x) = 9$ while training the three RLs and designing the backstepping using  $\gamma = 9$.}
	\label{fig:20}
\end{figure}

To evaluate robustness under model mismatch, we simulate a system with $\gamma = 8.5$ using the three LR\_based controllers trained at $\gamma =9$ and the backstepping controller designed for $\gamma = 9$. The results in Fig. \ref{fig:20} demonstrate that NOSAC\_training achieves superior robustness, outperforming the SAC and NOSAC controllers in overshoot, convergence speed, and steady-state error.

\section{Conclusion}
This paper proposes a novel neural-operator-embedded SAC architecture for PDE control, integrating rigorous control knowledge via a DeepONet approximating the backstepping controller. Applied to 1D hyperbolic and parabolic unstable PDEs, it is evaluated against backstepping control, standard SAC, and SAC with non-pretrained DeepONet in terms of overshoot, convergence rate, and steady-state error. Simulations show superior performance. Robustness to model mismatch (via modified system coefficients) is validated, with pre-trained DeepONet enabling adaptation to coefficient variations and generating appropriate control signals, ensuring better performance. Future work will explore integrating safety control with this learning-based approach.

\bibliographystyle{IEEEtran}
\bibliography{reference}

\end{document}